\documentclass[11 pt]{amsart}

\usepackage{amsmath}
\usepackage{amssymb}
\usepackage{amsfonts}
\usepackage{amsthm}
\usepackage{enumerate}
\usepackage[mathscr]{eucal}

\newtheorem{thm}{Theorem}[section]
\newtheorem{lem}[thm]{Lemma}
\newtheorem{sublem}[thm]{Sublemma}
\newtheorem{corollary}[thm]{Corollary}
\newtheorem{prop}[thm]{Proposition}

\numberwithin{equation}{section}

\begin{document}

\theoremstyle{definition}
\newtheorem{definition}{Definition}[section]
\newtheorem{remark}{Remark}[section]
\newtheorem{problem}{Problem}[section]

\numberwithin{equation}{section}

\def\intslash{\rlap{\kern  .32em $\mspace {.5mu}\backslash$ }\int}
\def\qsl{{\rlap{\kern  .32em $\mspace {.5mu}\backslash$ }\int_{Q_x}}}
\def\Re{\operatorname{Re\,}}
\def\Im{\operatorname{Im\,}}
\def\mx{{\max}}
\def\mn{{\min}}
\def\vth{\vartheta}

\def\rn{\rr^{n}}
\def\rr{\mathbb R}
\def\R{\mathbb R}
\def\Q{\mathcal Q}
\def\N{\mathbb N}
\def\complex{{\mathbb C}}
\def\norm#1{{ \left|  #1 \right| }}
\def\Norm#1{{ \left\|  #1 \right\| }}
\def\set#1{{ \left\{ #1 \right\} }}
\def\floor#1{{\lfloor #1 \rfloor }}
\def\emph#1{{\it #1 }}
\def\diam{{\text{\rm diam}}}
\def\dyad{{\text{\rm dyad}}}
\def\osc{{\text{\rm osc}}}
\def\ffB{\mathcal B}
\def\itemize#1{\item"{#1}"}
\def\seq{\subseteq}
\def\Id{\text{\sl Id}}

\def\Ga{\Gamma}
\def\ga{\gamma}
\def\Th{\Theta}

\def\prd{{\text{\it prod}}}
\def\parab{{\text{\it parabolic}}}

\def\eg{{\it e.g. }}
\def\cf{{\it cf}}
\def\Rn{{\mathbb R^n}}
\def\Rd{{\mathbb R^d}}
\def\sgn{{\text{\rm sign }}}
\def\rank{{\text{\rm rank }}}
\def\corank{{\text{\rm corank }}}
\def\coker{{\text{\rm Coker }}}
\def\loc{{\text{\rm loc}}}
\def\spec{{\text{\rm spec}}}

\def\comp{{\text{\rm comp}}}

\def\Coi{{C^\infty_0}}
\def\dist{{\text{\rm dist}}}
\def\diag{{\text{\rm diag}}}
\def\supp{{\text{\rm supp }}}
\def\rad{{\text{\rm rad}}}
\def\Lip{{\text{\rm Lip}}}
\def\trfrac#1#2{\tfrac{#2}{#1}}
\def\inn#1#2{\langle#1,#2\rangle}
\def\biginn#1#2{\big\langle#1,#2\big\rangle}
\def\rta{\rightarrow}
\def\lta{\leftarrow}
\def\noi{\noindent}
\def\lcontr{\rfloor}
\def\lco#1#2{{#1}\lcontr{#2}}
\def\lcoi#1#2{\imath({#1}){#2}}
\def\rco#1#2{{#1}\rcontr{#2}}
\def\bin#1#2{{\pmatrix {#1}\\{#2}\endpmatrix}}
\def\meas{{\text{\rm meas}}}

\def\card{\text{\rm card}}
\def\lc{\lesssim}
\def\gc{\gtrsim}
\def\pv{\text{\rm p.v.}}

\def\alp{\alpha}             \def\Alp{\Alpha}
\def\bet{\beta}
\def\gam{\gamma}             \def\Gam{\Gamma}
\def\del{\delta}             \def\Del{\Delta}
\def\eps{\varepsilon}
\def\ep{\epsilon}
\def\zet{\zeta}
\def\tet{\theta}             \def\Tet{\Theta}
\def\iot{\iota}
\def\kap{\kappa}
\def\ka{\kappa}
\def\lam{\lambda}            \def\Lam{\Lambda}
\def\la{\lambda}             \def\La{\Lambda}
\def\sig{\sigma}             \def\Sig{\Sigma}
\def\si{\sigma}              \def\Si{\Sigma}
\def\vphi{\varphi}
\def\ome{\omega}             \def\Ome{\Omega}
\def\om{\omega}              \def\Om{\Omega}

\def\fA{{\mathfrak {A}}}
\def\fB{{\mathfrak {B}}}
\def\fC{{\mathfrak {C}}}
\def\fD{{\mathfrak {D}}}
\def\fE{{\mathfrak {E}}}
\def\fF{{\mathfrak {F}}}
\def\fG{{\mathfrak {G}}}
\def\fH{{\mathfrak {H}}}
\def\fI{{\mathfrak {I}}}
\def\fJ{{\mathfrak {J}}}
\def\fK{{\mathfrak {K}}}
\def\fL{{\mathfrak {L}}}
\def\fM{{\mathfrak {M}}}
\def\fN{{\mathfrak {N}}}
\def\fO{{\mathfrak {O}}}
\def\fP{{\mathfrak {P}}}
\def\fQ{{\mathfrak {Q}}}
\def\fR{{\mathfrak {R}}}
\def\fS{{\mathfrak {S}}}
\def\fT{{\mathfrak {T}}}
\def\fU{{\mathfrak {U}}}
\def\fV{{\mathfrak {V}}}
\def\fW{{\mathfrak {W}}}
\def\fX{{\mathfrak {X}}}
\def\fY{{\mathfrak {Y}}}
\def\fZ{{\mathfrak {Z}}}

\def\fa{{\mathfrak {a}}}
\def\fb{{\mathfrak {b}}}
\def\fc{{\mathfrak {c}}}
\def\fd{{\mathfrak {d}}}
\def\fe{{\mathfrak {e}}}
\def\ff{{\mathfrak {f}}}
\def\fg{{\mathfrak {g}}}
\def\fh{{\mathfrak {h}}}
\def\fj{{\mathfrak {j}}}
\def\fk{{\mathfrak {k}}}
\def\fl{{\mathfrak {l}}}
\def\fm{{\mathfrak {m}}}
\def\fn{{\mathfrak {n}}}
\def\fo{{\mathfrak {o}}}
\def\fp{{\mathfrak {p}}}
\def\fq{{\mathfrak {q}}}
\def\fr{{\mathfrak {r}}}
\def\fs{{\mathfrak {s}}}
\def\ft{{\mathfrak {t}}}
\def\fu{{\mathfrak {u}}}
\def\fv{{\mathfrak {v}}}
\def\fw{{\mathfrak {w}}}
\def\fx{{\mathfrak {x}}}
\def\fy{{\mathfrak {y}}}
\def\fz{{\mathfrak {z}}}

\def\bbA{{\mathbb {A}}}
\def\bbB{{\mathbb {B}}}
\def\bbC{{\mathbb {C}}}
\def\bbD{{\mathbb {D}}}
\def\bbE{{\mathbb {E}}}
\def\bbF{{\mathbb {F}}}
\def\bbG{{\mathbb {G}}}
\def\bbH{{\mathbb {H}}}
\def\bbI{{\mathbb {I}}}
\def\bbJ{{\mathbb {J}}}
\def\bbK{{\mathbb {K}}}
\def\bbL{{\mathbb {L}}}
\def\bbM{{\mathbb {M}}}
\def\bbN{{\mathbb {N}}}
\def\bbO{{\mathbb {O}}}
\def\bbP{{\mathbb {P}}}
\def\bbQ{{\mathbb {Q}}}
\def\bbR{{\mathbb {R}}}
\def\bbS{{\mathbb {S}}}
\def\bbT{{\mathbb {T}}}
\def\bbU{{\mathbb {U}}}
\def\bbV{{\mathbb {V}}}
\def\bbW{{\mathbb {W}}}
\def\bbX{{\mathbb {X}}}
\def\bbY{{\mathbb {Y}}}
\def\bbZ{{\mathbb {Z}}}

\def\cA{{\mathcal {A}}}
\def\cB{{\mathcal {B}}}
\def\cC{{\mathcal {C}}}
\def\cD{{\mathcal {D}}}
\def\cE{{\mathcal {E}}}
\def\cF{{\mathcal {F}}}
\def\cG{{\mathcal {G}}}
\def\cH{{\mathcal {H}}}
\def\cI{{\mathcal {I}}}
\def\cJ{{\mathcal {J}}}
\def\cK{{\mathcal {K}}}
\def\cL{{\mathcal {L}}}
\def\cM{{\mathcal {M}}}
\def\cN{{\mathcal {N}}}
\def\cO{{\mathcal {O}}}
\def\cP{{\mathcal {P}}}
\def\cQ{{\mathcal {Q}}}
\def\cR{{\mathcal {R}}}
\def\cS{{\mathcal {S}}}
\def\cT{{\mathcal {T}}}
\def\cU{{\mathcal {U}}}
\def\cV{{\mathcal {V}}}
\def\cW{{\mathcal {W}}}
\def\cX{{\mathcal {X}}}
\def\cY{{\mathcal {Y}}}
\def\cZ{{\mathcal {Z}}}

\def\tA{{\widetilde{A}}}
\def\tB{{\widetilde{B}}}
\def\tC{{\widetilde{C}}}
\def\tD{{\widetilde{D}}}
\def\tE{{\widetilde{E}}}
\def\tF{{\widetilde{F}}}
\def\tG{{\widetilde{G}}}
\def\tH{{\widetilde{H}}}
\def\tI{{\widetilde{I}}}
\def\tJ{{\widetilde{J}}}
\def\tK{{\widetilde{K}}}
\def\tL{{\widetilde{L}}}
\def\tM{{\widetilde{M}}}
\def\tN{{\widetilde{N}}}
\def\tO{{\widetilde{O}}}
\def\tP{{\widetilde{P}}}
\def\tQ{{\widetilde{Q}}}
\def\tR{{\widetilde{R}}}
\def\tS{{\widetilde{S}}}
\def\tT{{\widetilde{T}}}
\def\tU{{\widetilde{U}}}
\def\tV{{\widetilde{V}}}
\def\tW{{\widetilde{W}}}
\def\tX{{\widetilde{X}}}
\def\tY{{\widetilde{Y}}}
\def\tZ{{\widetilde{Z}}}

\def\tcA{{\widetilde{\mathcal {A}}}}
\def\tcB{{\widetilde{\mathcal {B}}}}
\def\tcC{{\widetilde{\mathcal {C}}}}
\def\tcD{{\widetilde{\mathcal {D}}}}
\def\tcE{{\widetilde{\mathcal {E}}}}
\def\tcF{{\widetilde{\mathcal {F}}}}
\def\tcG{{\widetilde{\mathcal {G}}}}
\def\tcH{{\widetilde{\mathcal {H}}}}
\def\tcI{{\widetilde{\mathcal {I}}}}
\def\tcJ{{\widetilde{\mathcal {J}}}}
\def\tcK{{\widetilde{\mathcal {K}}}}
\def\tcL{{\widetilde{\mathcal {L}}}}
\def\tcM{{\widetilde{\mathcal {M}}}}
\def\tcN{{\widetilde{\mathcal {N}}}}
\def\tcO{{\widetilde{\mathcal {O}}}}
\def\tcP{{\widetilde{\mathcal {P}}}}
\def\tcQ{{\widetilde{\mathcal {Q}}}}
\def\tcR{{\widetilde{\mathcal {R}}}}
\def\tcS{{\widetilde{\mathcal {S}}}}
\def\tcT{{\widetilde{\mathcal {T}}}}
\def\tcU{{\widetilde{\mathcal {U}}}}
\def\tcV{{\widetilde{\mathcal {V}}}}
\def\tcW{{\widetilde{\mathcal {W}}}}
\def\tcX{{\widetilde{\mathcal {X}}}}
\def\tcY{{\widetilde{\mathcal {Y}}}}
\def\tcZ{{\widetilde{\mathcal {Z}}}}

\def\tfA{{\widetilde{\mathfrak {A}}}}
\def\tfB{{\widetilde{\mathfrak {B}}}}
\def\tfC{{\widetilde{\mathfrak {C}}}}
\def\tfD{{\widetilde{\mathfrak {D}}}}
\def\tfE{{\widetilde{\mathfrak {E}}}}
\def\tfF{{\widetilde{\mathfrak {F}}}}
\def\tfG{{\widetilde{\mathfrak {G}}}}
\def\tfH{{\widetilde{\mathfrak {H}}}}
\def\tfI{{\widetilde{\mathfrak {I}}}}
\def\tfJ{{\widetilde{\mathfrak {J}}}}
\def\tfK{{\widetilde{\mathfrak {K}}}}
\def\tfL{{\widetilde{\mathfrak {L}}}}
\def\tfM{{\widetilde{\mathfrak {M}}}}
\def\tfN{{\widetilde{\mathfrak {N}}}}
\def\tfO{{\widetilde{\mathfrak {O}}}}
\def\tfP{{\widetilde{\mathfrak {P}}}}
\def\tfQ{{\widetilde{\mathfrak {Q}}}}
\def\tfR{{\widetilde{\mathfrak {R}}}}
\def\tfS{{\widetilde{\mathfrak {S}}}}
\def\tfT{{\widetilde{\mathfrak {T}}}}
\def\tfU{{\widetilde{\mathfrak {U}}}}
\def\tfV{{\widetilde{\mathfrak {V}}}}
\def\tfW{{\widetilde{\mathfrak {W}}}}
\def\tfX{{\widetilde{\mathfrak {X}}}}
\def\tfY{{\widetilde{\mathfrak {Y}}}}
\def\tfZ{{\widetilde{\mathfrak {Z}}}}

\def\Atil{{\widetilde A}}          \def\atil{{\tilde a}}
\def\Btil{{\widetilde B}}          \def\btil{{\tilde b}}
\def\Ctil{{\widetilde C}}          \def\ctil{{\tilde c}}
\def\Dtil{{\widetilde D}}          \def\dtil{{\tilde d}}
\def\Etil{{\widetilde E}}          \def\etil{{\tilde e}}
\def\Ftil{{\widetilde F}}          \def\ftil{{\tilde f}}
\def\Gtil{{\widetilde G}}          \def\gtil{{\tilde g}}
\def\Htil{{\widetilde H}}          \def\htil{{\tilde h}}
\def\Itil{{\widetilde I}}          \def\itil{{\tilde i}}
\def\Jtil{{\widetilde J}}          \def\jtil{{\tilde j}}
\def\Ktil{{\widetilde K}}          \def\ktil{{\tilde k}}
\def\Ltil{{\widetilde L}}          \def\ltil{{\tilde l}}
\def\Mtil{{\widetilde M}}          \def\mtil{{\tilde m}}
\def\Ntil{{\widetilde N}}          \def\ntil{{\tilde n}}
\def\Otil{{\widetilde O}}          \def\otil{{\tilde o}}
\def\Ptil{{\widetilde P}}          \def\ptil{{\tilde p}}
\def\Qtil{{\widetilde Q}}          \def\qtil{{\tilde q}}
\def\Rtil{{\widetilde R}}          \def\rtil{{\tilde r}}
\def\Stil{{\widetilde S}}          \def\stil{{\tilde s}}
\def\Ttil{{\widetilde T}}          \def\ttil{{\tilde t}}
\def\Util{{\widetilde U}}          \def\util{{\tilde u}}
\def\Vtil{{\widetilde V}}          \def\vtil{{\tilde v}}
\def\Wtil{{\widetilde W}}          \def\wtil{{\tilde w}}
\def\Xtil{{\widetilde X}}          \def\xtil{{\tilde x}}
\def\Ytil{{\widetilde Y}}          \def\ytil{{\tilde y}}
\def\Ztil{{\widetilde Z}}          \def\ztil{{\tilde z}}


\def\ahat{{\hat a}}          \def\Ahat{{\widehat A}}
\def\bhat{{\hat b}}          \def\Bhat{{\widehat B}}
\def\chat{{\hat c}}          \def\Chat{{\widehat C}}
\def\dhat{{\hat d}}          \def\Dhat{{\widehat D}}
\def\ehat{{\hat e}}          \def\Ehat{{\widehat E}}
\def\fhat{{\hat f}}          \def\Fhat{{\widehat F}}
\def\ghat{{\hat g}}          \def\Ghat{{\widehat G}}
\def\hhat{{\hat h}}          \def\Hhat{{\widehat H}}
\def\ihat{{\hat i}}          \def\Ihat{{\widehat I}}
\def\jhat{{\hat j}}          \def\Jhat{{\widehat J}}
\def\khat{{\hat k}}          \def\Khat{{\widehat K}}
\def\lhat{{\hat l}}          \def\Lhat{{\widehat L}}
\def\mhat{{\hat m}}          \def\Mhat{{\widehat M}}
\def\nhat{{\hat n}}          \def\Nhat{{\widehat N}}
\def\ohat{{\hat o}}          \def\Ohat{{\widehat O}}
\def\phat{{\hat p}}          \def\Phat{{\widehat P}}
\def\qhat{{\hat q}}          \def\Qhat{{\widehat Q}}
\def\rhat{{\hat r}}          \def\Rhat{{\widehat R}}
\def\shat{{\hat s}}          \def\Shat{{\widehat S}}
\def\that{{\hat t}}          \def\That{{\widehat T}}
\def\uhat{{\hat u}}          \def\Uhat{{\widehat U}}
\def\vhat{{\hat v}}          \def\Vhat{{\widehat V}}
\def\what{{\hat w}}          \def\What{{\widehat W}}
\def\xhat{{\hat x}}          \def\Xhat{{\widehat X}}
\def\yhat{{\hat y}}          \def\Yhat{{\widehat Y}}
\def\zhat{{\hat z}}          \def\Zhat{{\widehat Z}}

\def\tg{{\widetilde g}}
\def\tpsi{{\widetilde {\psi}}}

\author{ Loukas Grafakos,  Petr Honz\'\i k,
 Andreas Seeger}
\address{
        L. Grafakos\\
        Department of Mathematics\\
        University of Missouri \\
        Co\-lum\-bia,  MO 65211, USA}
\email{loukas@math.missouri.edu}

\address{P. Honz\'\i{}k\\
        Department of Mathematics\\
        University of Missouri \\
        Columbia, MO 65211, USA}
\email{honzikp@math.missouri.edu}

\address{
        A. Seeger\\
        Department of Mathematics\\
        University of Wisconsin \\
        Madison, WI 53706, USA}
\email{seeger@math.wisc.edu}
\thanks{Grafakos and Seeger
were supported in part by  NSF grants.
Honz\'\i k was supported by 201/03/0931 Grant Agency of the
Czech Republic}

\title[On maximal functions for  Mikhlin-H\"ormander multipliers]
{On Maximal functions for  \\ Mikhlin-H\"ormander multipliers}

\date{October 26, 2004}

\begin{abstract}
Given
Mikhlin-H\"ormander multipliers $m_i$, $i=1,\dots, N$, with uniform estimates
we prove an  optimal $\sqrt{\log(N+1)}$ bound in $L^p$  for the maximal function
$\sup_i|\cF^{-1}[m_i\widehat f]|$ and related bounds for maximal functions
generated by dilations. These improve results in \cite{CGHS}.
\end{abstract}

\maketitle

\section{Introduction}\label{intro}

Given a symbol $m$ satisfying
\begin{equation}\label{symbols} |\partial^\alpha  m(\xi)|\le C_\alpha|\xi|^{-\alpha}
\end{equation}
for all multiindices $\alpha$,
then  by classical Calder\'on-Zygmund theory  the operator
$f\mapsto \cF^{-1}[m \widehat f]$ defines an $L^p$ bounded operator.
We study  two types of maximal operators associated to such symbols.

First we consider
 $N$ multipliers
$m_1,\dots, m_N$
  satisfying uniformly  the
conditions \eqref{symbols} and  ask for bounds
\begin{equation}\label{Nmultipliers}
\big\|\sup_{1\le i\le N}|\cF^{-1}[m_i \widehat f]|\big\|_p\le A(N) \|f\|_p,
\end{equation}
for all $f\in \cS$.

Secondly we form two  maximal functions generated by dilations of a
single multiplier,
\begin{align}
\cM_m^{\dyad} f(x)&= \sup_{k\in \bbZ}
|\cF^{-1}[m(2^k\cdot) \widehat f]|
\label{maxdyad}
\\
\cM_m f(x)&= \sup_{t>0}|\cF^{-1}[m(t\cdot) \widehat f]|
\label{maxcont}
\end{align}
 and ask under what additional  conditions  on $m$  these
define bounded ope\-rators on $L^p$.

Concerning  \eqref{maxdyad}, \eqref{maxcont}
a counterexample in \cite{CGHS} shows that
in general additional conditions on $m$ are needed for the maximal
 inequality to hold; moreover positive results
were shown using rather  weak decay assumptions on $m$.
The counterexample also  shows that the optimal uniform bound  in
\eqref{Nmultipliers} satisfies
\begin{equation}\label{lowerbound}
A(N)\ge c \sqrt{\log(N+1)}.
\end{equation}
The extrapolation argument in \cite{CGHS} only gives the upper bound
$A(N)=O(\log (N+1))$ and the main purpose of this paper is to close this gap
and to show that the upper bound is indeed $O(\sqrt{\log (N+1)})$.

We will formulate our theorems with
minimal smoothness assumptions that will be described now.

Let $\phi\in C^\infty_0(\bbR^d)$ be supported in
$\{\xi:1/2<|\xi|<2\}$ so that
$$\sum_{k\in \bbZ} \phi(2^{-k}\xi)=1$$ for all
$\xi\in \bbR^d\setminus\{0\}$.
Let $\eta_0\in C^\infty_c(\bbR^d)$ so that $\eta_0$ is even,
$\eta_0(x)=1$ for $|x|\le 1/2$ and $\eta_0$ is supported where
$|x|\le 1$. For $\ell>0$ let $\eta_\ell(x)=
\eta_0(2^{-\ell}(x))-\eta_0(2^{-\ell+1}x)$ and define
$$H_{k,\ell}[m](x)=\eta_\ell(x) \cF^{-1}[ \phi m(2^k\cdot)](x).$$
In what follows we set
$$\|m\|_{Y(q,\alpha)}:= \sup_{k\in \Bbb Z}\sum_{\ell\ge 0}
2^{\ell \alpha}\|H_{k,\ell}[m]\|_{L^q}.$$

Using the Hausdorff-Young inequality one gets
\begin{equation}\label{locbesov}
\|m\|_{Y(r',\alpha)}\lc \sup_{k\in\bbZ}
\|\phi m(2^k\cdot)\|_{B^{r}_{\alpha,1}},
\qquad \text{ if $1\le r\le 2$}
\end{equation}
where  $B^r_{\alpha,1}$
is the usual Besov space ({\it cf.} Lemma \ref{sbl4} below).
Thus if $m$ belongs to  $Y(2,d/2)$ then
$m$ is a Fourier multiplier of
$L^p(\Bbb R^d)$, for $1<p<\infty$ (this follows from a slight
modification of  Stein's approach in
\cite{Stein}, ch. IV.3, see also \cite{Sgarchive} for a related endpoint bound).


\begin{thm}\label{Thm1}
Suppose that $1\le r<2$ and suppose
that the multipliers
  $m_i$, $i=1,\dots, N$ satisfy the condition
\begin{equation}\label{assumpthm1}
\sup_i  \|m_i\|_{Y(r', d/r)}\le B<\infty.
\end{equation}
Then for $r< p<\infty$
$$\big\|\sup_{i=1,\dots,N}
\big |\cF^{-1}[m_i \widehat f]\big|\,
\big\|_p\le C_{p,r} B \sqrt{\log(N+1)} \|f\|_p.$$
\end{thm}


%

By \eqref{locbesov} we immediately get

\begin{corollary}\label{cor1}
Suppose that $1<r<2$,  and
\begin{equation}\label{besovassumption}
\sup_{1\le i\le N}\sup_{t>0} \|\phi m_i(t\cdot)\|_{B^r_{d/r,1}}\le A.
\end{equation}
Then for $r<p<\infty$
$$\big\|\sup_{i=1,\dots,N} \big |\cF^{-1}[m_i \widehat f]\big|\,\big\|_p
\le C_{p,r} A \sqrt{\log(N+1)} \|f\|_p.$$
\end{corollary}

\noi{\it Remark.} If one uses $Y(\infty, d+\eps)$  in \eqref{assumpthm1}
or $B^1_{d+\eps,1}$ in \eqref{besovassumption} one can use
Calder\'on-Zygmund theory (see \cite{DT}, \cite{CGHS})
to prove the $H^1-L^1$ boundedness
and the  weak type $(1,1)$ inequality, both with
constant $O(\sqrt{\log(N+1)})$.

%
%

\

Our second result is concerned with the operators
$\cM_m^\dyad$, $\cM_m$
generated by dilations.

\begin{thm}\label{Thm2}
Suppose $1<p<\infty$,
 $q=\min \{p,2\}$.

(i)
 Suppose that
\begin{equation}
\label{sob-Lq}
 \|\phi m(2^k \cdot)\|_{L^q_\alpha}\le \omega(k), \quad k\in \bbZ,
\end{equation}
holds for $\alpha>d/q$
and suppose that  the nonincreasing rearrangement $\omega^*$ satisfies
\begin{equation}
\label{omega}
\omega^*(0)+\sum_{l=2}^\infty\frac{\om^*(l)}{l \sqrt{\log l}}
  <\infty.
\end{equation}
Then $\cM_m^\dyad$
is bounded on $L^p(\bbR^d)$.

(ii)  Suppose that \eqref{omega} holds and
\eqref{sob-Lq} holds for $\alpha> d/p+1/p'$ if $1<p\le 2$ or for
$\alpha>d/2+1/p$ if $p>2$.
Then $\cM_m$
is bounded on $L^p(\bbR^d)$.

If
\eqref{sob-Lq}, \eqref{omega} are satisfied with $q=1$, $\alpha>d$ then
$\cM_m$ is of weak type $(1,1)$, and $\cM_m$ maps $H^1$ to $L^1$.
\end{thm}

This improves the earlier result in \cite{CGHS} where the conclusion
is obtained under the assumption
$\sum_{l=2}^\infty {\om^*(l)}/{l}<\infty$, however somewhat weaker
smoothness assumptions  were made in \cite{CGHS}.

In \S\ref{model} we shall discuss  model cases  for
Rademacher expansions. In \S\ref{sect2} we shall give the outline of
the proof of Theorem
\ref{Thm1} which is based on the $\exp(L^2)$ estimate  by Chang-Wilson-Wolff \cite{CWW},  for functions with bounded Littlewood-Paley square-function.
The proof of a critical pointwise inequality is given in \S\ref{sect3}.
The proof of Theorem \ref{Thm2} is sketched in \S\ref{PfThm2}.
Some open problems are mentioned in
 \S\ref{open}.

\medskip

{\it Acknowledgement:} The second named author would like to thank Lubo\v s  Pick for a
helpful conversation concerning convolution inequalities in rearrangement
 invariant function spaces.

\section{Dyadic  model cases for Rademacher expansions} \label{model}
Before we discuss
 the proof of Theorem \ref{Thm1} we give a simple result on
expansions for Rademacher functions $r_j$ on $[0,1]$  which motivated the
proof.

\begin{prop}\label{P1}
Let $a^i\in \ell^2$.
and let
$$ F_i(s)=\sum_j a_j^i r_j(s), \quad s\in [0,1].
$$
Then
$$
\big\|\sup_{i<N}|F_i|\big\|_{L^2[0,1]}\lc \sup\|a^i\|_{\ell^2}
\sqrt{\log (N+1)}.
$$
\end{prop}

\begin{proof}
We use  the well known estimate for the distribution function of
the Rademacher expansions (\cite{Stein}, p. 277),
\begin{equation}\label{Rademacher}
\meas\big(\{s\in [0,1]:|F_i(s)|>\la\}\big)\leq
2\exp\big(-\tfrac{\la^2}{4\|a^i\|_{\ell^2}^2}\big)
\end{equation}
Set $u_N=(4\log (N+1))^{1/2}\sup_{1\le i\le N}\|a^i\|_{\ell^2}$.
Then
\begin{align*}
\big\|&\sup_{i=1,\dots,N}|F_i|\big\|_2^2 \le
u_N^2+2\sum_{i=1}^N\int_{u_N}^{\infty}\la\,
\meas\big(\{s:|F_i(s)|>\la\}\big) d\la
\\&\leq
u_N^2+4\sum_{i=1}^N\int_{u_N}^{\infty}\la e^{-\la^2/(4\|a^i\|^2_{\ell^2})} d
\alpha\leq u_N^2 + 4\sup_{i=1,\dots, N}\|a^i\|_{\ell^2}^2 Ne^{-u_N^2/4}
\end{align*}
which is bounded by $(1+4 \log(N+1)) \sup_i\|a^i\|_{\ell^2}^2$. The
claim follows.
\end{proof}

There is a multiplier interpretation to this inequality. One
 can work with a single   function $f=\sum a_j r_j$ and  a family
of bounded sequences (or multipliers) $\{b^i\}$ and
one forms $F_i(s)=
\sum_{j} b^i_j a_j r_j(s)$.
 The norm then grows as a square root of the
logarithm of the  number of  multipliers; i.e. we have
\begin{corollary}
$$
\Big\| \sup_{i=1,\dots, N}\big|\sum_j b^i_j a_j r_j
\big|\Big\|_{L^2([0,1])}\lc
\sup_i\|b^i\|_\infty \sqrt{\log(N+1)}
\Big\| \sum_j  a_j r_j\Big\|_{L^2([0,1])}.
$$
\end{corollary}

We shall now consider a dyadic model case for the maximal operators generated by dilations.

\begin{prop}\label{Rademacher-dilation}
Consider a sequence $b=\{b_i\}_{i\in \bbZ}$ which satisfies
 $$b^*(l)\leq \frac{A}{(\log (l+2))^{1/2}}.$$
Then for any sequence $a=\{a_n\}_{n=1}^\infty$ we have
$$\Big\|\sup_{k\in \bbZ} \big|\sum_{j=0}^{\infty} b_{j-k}a_j r_j\big|
\Big\|_2\leq C A\|a\|_2.
$$
\end{prop}

\begin{proof} We may assume that both $a$ and $b$ are real valued sequences.
Let
$$H_k(s)=\sum_{j=1}^{\infty} b_{j-k}a_j r_j(s).
$$
Then by orthogonality of the Rademacher functions
$$\|H_k\|_2^2=\sum_{j=1}^{\infty} [b_{j-k}a_j]^2.$$
We shall use a result of Calder\' on~\cite{Calderon}
which states  that if some linear operator is bounded on $L^1(\mu)$ and
on $L^{\infty}(\mu)$ on a space with $\sigma$-finite measure $\mu$,
then it is bounded on all rearrangement invariant function spaces
on that space.
In our case the intermediate space is the Orlicz space ${\exp}\, \ell$,
which coincides with
the space of all sequences $\gamma=\{\gamma_j\}_{j\in \Bbb Z}$  that satisfy
the condition \begin{equation}\label{explnorm}
\gamma^*(l)\leq \frac{C}{\log (l+2)}, \quad l\ge 0,
\end{equation} and the best constant in \ref{explnorm} is equivalent
to the norm in
$\exp(\ell)$. We apply Calder\'on's result to the operator $T$ defined by
$$
[T\gamma]_k= \sum_{j=1}^\infty \gamma_{j-k} a_j^2
$$
and get
\begin{equation*}
\sup_{l\ge 0} \log(l+2) (T\gamma)^*(l)\le C
\big\|\{a_n^2\}\big
\|_{\ell^1}
\sup_{l\ge 0} \log(l+2) \gamma^*(l).
\end{equation*}
Let $c_k= \|H_k\|_2
\equiv ([T (b^2)]_k)^{1/2}$
 where $b^2$ stands for the sequence $\{b_j^2\}$; then
by our bound for $T\gamma$  and the  assumption on $b$  it follows that
\begin{equation}\label{cjbound}
c^*(l)\le C_1 A \|a\|_{\ell^2}\big(\log(2+l)\big)^{-1/2}.
\end{equation}

We can proceed with the proof as in Proposition~\ref{P1}, using
again \eqref{Rademacher}, {\it i.e.}
$$\meas(\{s\in [0,1]:|H_k(s)|>\alpha\}) \leq 2e^{-\alpha^2/4c_k^2}.$$
Then we obtain for $u>0$
\begin{align}
\big\|\sup_k|H_k|\big\|_2&\le u^2 + 4\sum_k
\int_u^\infty \alpha e^{-\alpha^2/4c_k^2}
\notag
\\
&\le u^2 + 8\sum_{k} c_k^2 e^{-u^2/(4c_k^2)}\notag
\\&= u^2 + 8\sum_{l\ge 0}  (c^*(l))^2 e^{-u^2/4(c^*(l))^2}.\notag
\end{align}
We set the cutoff level to be $u=10 C_1 A\|a\|_2$
and obtain
\begin{align*}
&\|\sup_k|H_k|\|_2^2\le u^2 + C_1^2 A^2
\sum_{l\ge 0} (2+l)^{-5/2}
\lc A^2 \|a\|^2_2
\end{align*}
which is what we wanted to prove.
\end{proof}

\noi{\it Remark:} Since the $L^p$ norm of $\sum a_jr_j$ is equivalent to the
$\ell^2$ norm of $\{a_j\}$ one can also prove $L^p$
analogues of the two propositions,  for $0<p<\infty$.

\medskip

\section{Proof of Theorem \ref{Thm1}}\label{sect2}
To prove \eqref{Nmultipliers} we may assume that $\widehat f$ is
 compactly supported  in $\Bbb R^d\setminus \{0\}$ and thus we may assume that the multipliers $m_i$ are compactly
supported on a finite union of dyadic annuli. In view of the scale invariance of the assumptions we may assume without loss of generality that
\begin{equation}\label{vanish}
m_i(\xi)=0,\quad|\xi|\le 2^{N}, i=1,\dots, N.
\end{equation}

In the case of Fourier multipliers the  inequality \eqref{Rademacher}
will be replaced by a  ``good-$\lambda$ inequality'' involving
square-functions for martingales as proved by Chang, Wilson and Wolff
\cite{CWW}.
To fix notation let, for any $k\ge 0$,
$\fQ_k$ denote the family of dyadic cubes of sidelength $2^{-k}$; each $Q$ is
of the form $\prod_{i=1}^d [n_i 2^{-k}, (n_i+1) 2^{-k})$.
Denote by $\bbE_k$ the conditional expectation,
$$
\bbE_k f(x) =\sum_{Q\in\fQ_k} \chi_Q(x) \frac{1}{|Q|}\int_Q f(y) dy
$$
and by $\bbD_k$ the martingale differences,
$$
\bbD_k f(x)=\bbE_{k+1} f(x)-\bbE_{k} f(x).
$$
The square function for the dyadic martingale is defined by
$$
S(f)=\Big(\sum_{k\ge 0}|\bbD_k f(x)|^2\Big)^{1/2};
$$
one has the inequality $\|S(f)\|_p\le C_p\|f\|_p$ for $1<p<\infty$
(see \cite{BDG}, \cite{Bu1} for the general martingale case, and
for our special case  \cf.  also  Lemma \ref{lemma1} below).

The result from \cite{CWW} says that
there is a constant $c_d>0$ so that for all $\la>0$, $0<\eps<1$, one has

\begin{multline}\label{cww-inequality}
\meas\big(\big\{x:\sup_{k\ge 0} |\bbE_k g(x)-\bbE_0 g(x)|>2\lambda,\, S(g)< \epsilon\lambda\big\}\big))
\\ \leq
C{\rm exp}(-\frac{c_d}{\epsilon^2})
\meas\big(\big\{x:\sup_{k\ge 0} |\bbE_k g(x)|>\eps\lambda\big\}\big);
\end{multline}
see \cite{CWW} (Corollary 3.1 and a remark on page 236).
To use  \eqref{cww-inequality}
 we need a pointwise inequality for square functions
applied to convolution operators.

Choose a radial
Schwartz function $\psi$ which equals $1$ on the support of $\phi$ (defined in the introduction)  and is compactly supported in $\bbR^d\setminus \{0\}$, and define the Littlewood-Paley operator $L_k$ by
\begin{equation}\label{littlewoodpaley}\widehat {L_k f}(\xi)= \psi(2^{-k}\xi) m(\xi) \widehat f(\xi)
\end{equation}
Let $M$ be the Hardy-Littlewood maximal operator and define  the operator
$M_r$ by
$$M_r=(M(|f|^r))^{1/r}.$$ Denote by
$\fM= M\circ M\circ M$ the three-fold iteration of
the  maximal operator.
Now define
\begin{equation}\label{Grdefinition}
G_{r}(f)=
\Big(\sum_{k\in \bbZ}
\big(\fM[|L_k f|^r]\big)^{2/r}\Big)^{1/2}.
\end{equation}

{}From the Fefferman-Stein inequality for vector-valued maximal functions
\cite{FS},
\begin{equation} \label{FeffSt}
\|G_r(f)\|_p\le C_{p,r}\|f\|_p, \quad 1<r<2, \,r<p<\infty.
\end{equation}

\begin{lem}\label{lemma1}
Let $Tf=\mathcal F^{-1} [m\widehat f]$ and let $1< r\le \infty$.
 Then for $x\in \bbR^d$,
\begin{equation}
S(T f)(x) \le A_{r} \|m\|_{Y(r', d/r)} G_r(f)(x).
\end{equation}
\end{lem}

The proof will be given in \S\ref{sect3}.

We shall also need
\begin{lem}\label{errorlemma} Let $Tf=\cF^{-1}[m\widehat f]$ and
suppose that $m(\xi)=0$ for $|\xi|\le 2^{N}$. Then
\begin{equation}
|\bbE_0 Tf(x)|\le C2^{-N/r}
 C_r \|m\|_{Y(r', d/r)}  (\fM(|f|^r))^{1/r}.
\end{equation}
\end{lem}

\medskip

We  now give the proof of  Theorem \ref{Thm1}.
Let $T_i f=\cF^{-1}[m_i\widehat f]$.
We need to estimate
$$
\big\|\sup_{1\le i\le N}|T_i f|\big\|_p= \Big(p4^p\int_0^\infty
\lambda^{p-1}
\meas (\{x:\sup_i|T_i f(x)|>4\lambda\}) d\la\Big)^{1/p}.
$$
Now by Lemma \ref{lemma1} one gets the pointwise bound
\begin{equation}\label {ptw} S(T_i f)\le A_r B\, G_r(f).
\end{equation} We note that
$$\{x:\sup_{1\le i\le N}|T_i f(x)|>4\lambda\} \subset E_{\la,1}
\cup E_{\la,2}
\cup E_{\la,3}
$$
where with
\begin{equation}\label{epsilon}\eps_N:=
 \Big(\frac{c_d}{10\log(N+1)}\Big)^{1/2}\end{equation} we have set
\begin{align*}E_{\la,1}&=\{x: \sup_{1\le i\le N}
|T_i f(x)-\bbE_0 T_i f(x)|>2\la, G_r(f)(x)\le
\frac{\eps_N  \la}{A_r B}\},
\\
E_{\la,2}&=\{x:  G_r(f)(x)>
\frac{\eps_N  \la}{A_r B}\},
\\
E_{\la,3}&=\{x: \sup_{1\le i\le N}|\bbE_0 T_i f
(x)|>2\la\}.
\end{align*}
By \eqref{ptw},
\begin{equation}
E_{\la,1}\subset \bigcup_{i=1}^N
\{x: |T_i f(x)|>2\la, S(T_i f)\le \eps_N \la\},
\end{equation}
and thus using the good-$\lambda$ inequality \eqref{cww-inequality}
we obtain
\begin{align*}
\meas(E_{\la,1}) &\le \sum_{i=1}^N
\meas\big(\{x: |T_i f(x)-\bbE_0 T_if(x)|>2\la, S(T_i f)\le \eps_N \la\}\big)
\\
&\le \sum_{i=1}^N C
\exp(-\trfrac{\eps_N^2}{c_d})
\meas(\{x: \sup_k |\bbE_k (T_i f)|>\la\}).
\end{align*}
Hence
\begin{align}
&\Big(p\int_0^\infty \la^{p-1} \meas(E_{\la,1}) d\la\Big)^{1/p}
\notag
\\
&\lc \Big(\sum_{i=1}^N
 \exp(-\trfrac{\eps_N^2}{c_d})
\big\|\sup_k |\bbE_k (T_i f)|\big\|_p^p\Big)^{1/p}
\notag
\\
&\lc \Big(\sum_{i=1}^N
 \exp(-\trfrac{\eps_N^2}{c_d})\big\|T_i f\big\|_p^p\Big)^{1/p}
\notag
\\
&\lc B \big(N \exp(-\trfrac{\eps_N^2}{c_d})\big)^{1/p}\|f\|_p \lc B\|f\|_p
\label{integrateE1}
\end{align}
uniformly in $N$ (by our choice of $\eps_N$ in \eqref{epsilon}).

Next, by a change of variable,
\begin{align}
\Big(p\int_0^\infty \la^{p-1} \meas(E_{\la,2}) d\la\Big)^{1/p}
&= \frac{A_r B}{\eps_N} \big\|G_r(f)\big\|_p
\notag
\\&\lc B \sqrt{\log(N+1)}
\|f\|_p\label{integrateE2}
\end{align}

Finally, from Lemma \ref{errorlemma} and the Fefferman-Stein inequality
$$\meas(E_{\la,3}) \le
\sum_{i=1}^N
\meas\big(\{x: |\bbE_0 T_if(x)|>2\la\}\big)$$
and thus
\begin{align}
&\Big(p\int_0^\infty \la^{p-1} \meas(E_{\la,3}) d\la\Big)^{1/p}
= 2 \big\|\sup_{i=1,\dots,N}|\bbE_0(T_i f)|\big\|_p
\notag
\\
&\le 2\Big(\sum_{i=1}^N \big\|\bbE_0(T_i f)\big\|_p^p\Big)^{1/p}
\lc B N^{1/p} 2^{-N/r }\|f\|_p \lc B \|f\|_p.
\label{integrateE3}
\end{align}
The asserted inequality follows from
\eqref{integrateE1},
\eqref{integrateE2},
and \eqref{integrateE3}. \qed

\

For completeness we mention the well known  relation of the $Y(r',\alpha)$
conditions with Besov and Sobolev norms.

\begin{lem}\label{sbl4}
Let $1\le r\le 2$ and $\alpha>d/r$.
Then
\begin{equation*}
\begin{aligned}
\|m\|_{Y(r',d/r)}&\lc \sup_k\|\phi m(2^k\cdot)\|_{B^r_{d/r,1}}
\\&\lc
\sup_k\|\phi m(2^k\cdot)\|_{L^r_\alpha}
\lc \sup_k\|\phi m(2^k\cdot)\|_{L^2_\alpha}
\end{aligned}
\end{equation*}\end{lem}
\begin{proof}
By the Hausdorff-Young inequality and the definition of the Besov space
we have
$$ \sum_{\ell=0}^\infty 2^{\ell d/r}\|H_{k,\ell}\|_{r'}
\lc \sum_{\ell=0}^\infty 2^{\ell d/r}\|[\phi m(2^k\cdot) ]*\widehat{\eta_\ell}\|_{r}
\lc \|\phi m(2^k\cdot) \|_{B^r_{d/r,1}}.$$
By elementary  imbedding properties
$\|g\|_{B^r_{d/r,1}}\lc
 \|g\|_{L^r_\gamma}$
 if $\gamma>d/r$.
Finally
$\|\phi m(2^k\cdot) \|_{L^r_\gamma}
\lc C_r'  \|\phi m(2^k\cdot) \|_{L^2_\gamma}$,  if $1<r\le 2$.
In this last inequality we used that for $\chi\in C^\infty_c$ we have
$\|\chi g\|_{L^{r_0}_\gamma }
\lc \|g\|_{L^{r_1}_\gamma}$ for $r_0\le r_1$, $\gamma\ge 0$;
this is trivial for integers  $\gamma$ from H\"older's inequality
 and follows for all
$\gamma\ge 0$ by interpolation.
\end{proof}

\section{Proofs of Lemma \ref{lemma1}
and Lemma \ref{errorlemma}
}\label{sect3}

Choose a radial
Schwartz function $\beta$ with the property that $\widehat \beta$ is supported in $\{x:|x|\le 1/4\}$ so that $\beta(\xi)\neq 0$ in
$\{\xi:1/4\le|\xi|\le 4\}$ and $ \beta(0)=0$. Now choose a
function $\tpsi \in C^\infty_c$ so that $\tpsi(\xi)(\beta(\xi))^2=1$
for all $\xi \in \supp \phi$, here $\phi$ is as in the formulation of the theorem. Define operators $T_k$, $B_k$, $\tL_k$ by
\begin{align*}
\widehat {T_k f}(\xi)&= \phi(2^{-k}\xi) m(\xi) \widehat f(\xi)
\\
\widehat {B_k f}(\xi)&= \beta(2^{-k}\xi)  \widehat f(\xi)
\\
\widehat {\tL_k f}(\xi)&= \tpsi(2^{-k}\xi)  \widehat f(\xi).
\end{align*}
Then $T =\sum_k T_k = \sum_k B_k^2 \tL_k T_k L_k $
and we write
\begin{equation}\label{splitting}\Bbb D_k Tf=
\sum_{n\in \Bbb Z}\
( \Bbb D_k
 B_{k+n}) ( B_{k+n}
\tL_{k+n})  T_{k+n}L_{k+n}f.
\end{equation}

\begin{sublem}\label{sbl1}
\begin{equation}|B_k\tL_k f(x)|\lc Mf (x).
\end{equation}
\end{sublem}
\begin{proof} Immediate.\end{proof}

\begin{sublem}\label{sbl2}
For $s\ge 0$,
\begin{align}\label{large-s-first}
|\bbE_{k+1} B_{k+s} f(x)|+
|\bbE_{k} B_{k+s} f(x)|
&\lc 2^{-s/q'} M_qf(x)
\end{align}
and
\begin{equation}
\label{cancellation-s}
|\bbD_k B_{k-s} f(x)|\lc 2^{-s} Mf(x).
\end{equation}
\end{sublem}
\begin{proof}
We give the proof although the estimates are rather standard (for
similar calculations in other contexts see for  example \cite{Christ}, \cite{JKRW},
 \cite{GK}, \cite{JSW}).

For  \eqref{large-s-first}
first note  this inequality is trivial if $s$ is small and assume, say,
 $s\ge 10$.
For
$Q\in \fQ_k$, $s>0$ let
$b_s(Q)$  be the set of all
$x\in Q$  for which
the $\ell^\infty$ distance to the boundary of $Q$ is
$\le 2^{-k-s+1}$.

Fix a cube $Q_0\in \fQ_{k+1}$.
If $Q'$ is a dyadic subcube of sidelength $2^{-k-s+1}$  subcube which is not contained in $b_s(Q)$ then
$B_{k+s}[f\chi_{Q'}]$ is supported in  $Q_0$
and using the cancellation of $\mathcal F^{-1}[\beta]$
we see that $\bbE_{k+1} B_{k+s}[\chi_{Q'} g]=0$ for all $g$.
Let
$\cV_s(Q_0)$ be the union over all dyadic
cubes of sidelength $2^{-k-s+1}$ whose closures intersect
the boundary of $Q_0$. Then $$\bbE_{k+1} B_{k+s}[\chi_{Q_0} g]=
\bbE_{k+1} B_{k+s}[g \chi_{\cV_s(Q_0)}]$$  for all $g$.
In view of the support
properties of $\widehat \beta$ we note that
$B_{k+s}[g\chi_{\cV_s(Q_0)}]$ is also supported in $\cV_{s-1}(Q_0)$. Observe that this set has measure
$O(2^{-kd} 2^{-s})$.

It follows that for $x\in Q_0$
\begin{align*} |\bbE_{k+1}B_{k+s} f(x)|&\le 2^d |Q_0|^{-1}\int_{\cV_{s-1}(Q_0)}
|B_{k+s}[
\chi_{\cV_s(Q_0)} f](y) |dy
\\&\lc |Q_0|^{-1}\big(\int_{Q_0}
 |f(y) |^qdy\big)^{1/q} 2^{-(kd+s)/q'}
\\
&\lc 2^{-s/q'} \big(M(|f|^q)\big)^{1/q}
\end{align*}
By the same argument one obtains this bound also for
$|\bbE_k B_{k+s} f|$ and  thus
\eqref{large-s-first}
follows.

The inequality \eqref{cancellation-s}
$\bbD_k B_{k-s} f$ is a  simple consequence of the smoothness of the convolution kernel of $B_{k-s}$ and the cancellation properties of the
operator $\bbD_k=\bbE_{k+1}-\bbE_k$.
 \end{proof}

\begin{sublem}\label{sbl3}
Let $1<r<\infty$.
We have
\begin{equation}
|T_k f(x)|\le C\|m\|_{Y(r',d/r)}
  M_rf(x).
\end{equation}
\end{sublem}
\begin{proof}

We may decompose $T_k $ using the kernels $H_{k,l}$ and obtain
\begin{align*}
|T_k f(x)|&=\Big|\sum_{\ell=0}^\infty \int 2^{kd} H_{k,\ell}(2^ky) f(x-y) dy\Big|
\\&\le \sum_{\ell=0}^\infty
\Big(2^{kd} \int |H_{k,\ell}(2^k y) |^{r'} dy\Big)^{1/r'} \Big(2^{kd}\int_{|y|\le 2^{-k+\ell}}|f(x-y)|^r dy\Big)^{1/r}
\\&\le
 \sum_{\ell=0}^\infty 2^{\ell d/r}
\|H_{k,\ell}\|_{r'} \big(M(|f|^r)(x)\big)^{1/r}.\qed
\end{align*}


\medskip
\noindent{{\it Proof of Lemma \ref{lemma1}.}}
To estimate the terms in \eqref{splitting} we use
Sublemma \ref{sbl1} to bound
$B_{k+n}\tL_{k+n}$,
Sublemma \ref{sbl2}  to  bound
$\bbD_{k}B_{k+n}$ and
Sublemma \ref{sbl3}  to bound
$T_{k+n}$.
This yields that
\begin{multline*}|\Bbb D_k
 B_{k+n}^2 \tL_{k+n} T_{k+n} L_{k+n}f(x)|
\lc \|m\|_{Y(r',d/r)}\\ \times
\begin{cases}
2^{-n/q'} M_q\circ M\circ M_r
(L_{k+n} f)(x)&
\text{ if } n\ge 0
\\
2^{n}  M\circ M\circ M_r (L_{k+n} f)(x)
&\text{ if } n<0,
\end{cases}
\end{multline*}
and straightforward estimates imply the asserted bound.
\end{proof}

\medskip

\noindent{\it Proof of Lemma \ref{errorlemma}.}
We split $\bbE_0 Tf= \sum_{k\ge N-2}  \bbE_0
B_k^2 \tL_k T_k $, and
by
the sublemmas we get
$$
|\bbE_0 B_k^2 \tL_k T_k f(x)|\lc 2^{-k/r}  \|m\|_{Y(r',d/r)}
M_r\circ M\circ M_r(f) (x)
$$ which implies the assertion. \qed

\section
{Maximal functions generated by dilations}\label{PfThm2}
For the proof of Theorem \ref{Thm2} we use arguments in \cite{CGHS} and
applications of  Theorem \ref{Thm1}. Let us first  consider  the dyadic maximal operator $\cM_m^\dyad$.

Let $$\cI_j= \{k\in \bbZ: \omega^*(2^{2^j})<|\omega(k)|\le
\omega^*(2^{2^{j-1}})\}.$$
We split $m=\sum_j m_j$ where  $m_j$ is supported in the union of dyadic annuli
$\cup_{k\in \cI_j} \{\xi: 2^{k-1}<|\xi|<2^{k+1}\}$.

By Lemma 3.1 in \cite{CGHS} we can find a sequence of integers
$B=\{i\}$ so that
 for each $j$  the sets $b_i+\cI_j$ are pairwise disjoint, and
$\mathbb Z=\cup_{n=-4^{2^j+1}}^{4^{2^j+1}} (n+B)$.

\medskip

Let $T^j_k f=\cF^{-1}[ m_j(2^{k}\cdot)\widehat f]$.
We write
\begin{equation}\label{usingthebi}
\sup_k |T_k f|=
\sup_{|n|\le 4^{2^j+1}} \sup_{i\in \Bbb Z}|T_{b_i+n} f|
\end{equation}
and split the sup in $i$ according to
whether $i>0$, $i=0$, $i<0$.
We use the standard equivalence of the $L^p$ norm of expansions of Rademacher functions
 $\{r_i\}_{i=1}^\infty$ with the $\ell^2$ norm of the sequence of coefficients
(see \cite{Stein}, p. 276).

Then
\begin{align*}
\Big\|\sup_{|n|\le 4^{2^j+1}} \sup_{i>0}|T^j_{b_i+n} f|&\Big\|_p
\le\Big\|\sup_{|n|\le 4^{2^j+1}} \Big(\sum_{i>0}|T^j_{b_i+n} f|^2\Big)^{1/2}
\Big\|_p
\\
&\le C_p
\Big\|\sup_{|n|\le 4^{2^j+1}} \Big(\int_0^1\Big|\sum_{i=1}^\infty r_i(s)
 T^j_{b_i+n} f\Big|^p ds\Big)^{1/p} \Big\|_p
\\&\le C_p
\Big\|
\Big(\int_0^1 \sup_{|n|\le 4^{2^j+1}}
\Big|\sum_{i=1}^\infty r_i(s)
 T^j_{b_i+n} f\Big|^p ds\Big)^{1/p} \Big\|_p
\\
&= C_p
\Big(\int_0^1 \Big\|
 \sup_{|n|\le 4^{2^j}}
\Big|\sum_{i=1}^\infty r_i(s)
 T^j_{b_i+n} f\Big| \Big\|_p^p  ds\Big)^{1/p}
\end{align*}
which reduce matters for the dyadic maximal function
to an application of Theorem
\ref{Thm1}
(of course the terms above with $i\le 0$ are handled similarly).
Thus we obtain the estimate
$$\| M^\dyad_{m_j} \|_{L^p\to L^p} \lc
2^{j/2}\omega^*(2^{2^{j-1}}).$$

For the full maximal operator we use  standard decompositions by smoo\-thing
 out the rescaled dyadic pieces. We just sketch the argument.
Assume that $p\ge 2$ and that the
assumption  of Theorem \ref{Thm2}, (ii), with $\alpha>d/2+1/p$ holds. Then
one can  decompose $m_j=\sum_{l\ge 0} m_{j,l}$ where $m_{j,l}$
has essentially the same support property as $m_j$
(with slightly extended dyadic annuli) and where
$$\|\phi m_{j,l}(2^k\cdot)\|_{L^2_{\alpha-1/p}}
+
2^{-l}\|\phi  \,\inn{\xi}{\nabla} [m_{j,l}(2^k\cdot)]\|_{L^2_{\alpha-1/p}}
\lc \omega^*(2^{2^{j-1}}) 2^{-l/p}.
$$
One then uses a standard argument (see {\it e.g.} \cite{Stein3},  p. 499)
to see that
\begin{align*}
&\sup_{t>0}|\cF^{-1}[m_{j,l}(t\cdot) \widehat f]|
\,\le\,
C\sup_{k>0} |\cF^{-1}[m_{j,l}(2^k\cdot) \widehat f]| +
\\ &\, C
\Big(\int_{1}^2|\cF^{-1}[m_{j,l}(2^ku \cdot) \widehat f]|^pdu
\Big)^{\frac{1}{p'p}}
\Big(\int_{1}^2\big| (\partial/\partial u)\cF^{-1}[m_{j,l}(2^ku \cdot)
\widehat f]|^pdu
\Big)^{\frac{1}{p^2}}
\end{align*}
 and straightforward estimates reduce matters to the dyadic case
treated above. For the weak-type estimate (or the $H^1\to L^1$ estimate)
one has to combine this
argument with Calder\'on-Zygmund theory and the $L^p$ estimates for $1<p<2$ follow then by
   an analytic interpolation.
Similar arguments appear in
\cite{DT} and \cite{CGHS}; we omit the details.\qed

\section{Open problems}
\label{open}

Concerning Theorem \ref{Thm1} one can ask about $L^p$ boundedness for $p>2$
under merely the  assumption
 $m_i\in Y(p',\alpha)$,
 $\alpha >d/p$. Combining our present result with those in \cite{CGHS} one can show that if for some
 $2<r<\infty$
\begin{equation}
\label{rgreatertwo}\sup_{i}\|m_i\|_{Y(r',\alpha)}\le A, \quad   \alpha >d/r\end{equation}
then for $r\le p<\infty$
\begin{equation}\label{logbound}
\big\|\sup_{i=1,\dots,N} \big |\cF^{-1}[m_i \widehat f]\big|\,\big\|_p
\le C_{p,r,\alpha} A (\log(N+1))^{1/r'}
\|f\|_p.
\end{equation}
Indeed one can imbed the multipliers in analytic families
so that for $L^\infty\to BMO$ boundedness
 one has $Y(1+\eps_1, \eps_2)$ conditions and
the $O(\log (N+1))$
result  of \cite{CGHS} applies. For $p=2$ on has
the usual
$Y(2, d/2+\eps)$ conditions  and  Theorem \ref{Thm1} applies giving an
 $O((\log (N+1))^{1/2})$  bound.

\medskip

{\it Problem 1}:
Does \eqref{logbound} hold with an
$O(\sqrt{\log (N+1)})$ bound if we assuming \eqref{rgreatertwo} with
$r>2$ ?

\medskip

{\it Problem 2}: To which extent can one relax the smoothness conditions in
Theorems \ref{Thm1} and \ref{Thm2} to obtain $L^2$ bounds?
In particular what happens in Theorem \ref{Thm2} if one imposes
localized $L^2_\alpha$ conditions for $\alpha<d/2$,
 assuming again minimal decay assumptions on $\om^*$.

\medskip

Finally we discuss  possible optimal decay estimates
 for the maximal operators generated by dilations.
The hypothesis in Theorem \ref{Thm2} is equivalent with the  assumption
$$\{2^{j/2} \omega^*(2^{2^j})\}
\in \ell^1.$$
The counterexamples  in \cite{CGHS} leave open the possibility that the \
conclusion of Theorem \ref{Thm1} might hold under the weaker assumption
$\{2^{j/2} \omega^*(2^{2^j})\}\in \ell^\infty$, {\it i.e.}
\begin{equation}\label{log12decay}
 \omega^*(l)\le C \big(\log(2+l)\big)^{-1/2};
\end{equation}
this is in fact suggested by the dyadic model case in Proposition
\ref{Rademacher-dilation}.
The latter condition would be optimal  and leads us to formulate

\medskip

{\it Problem 3.}  Suppose $m$  is a symbol satisfying
\eqref{sob-Lq} for sufficiently large $\alpha$.
Does $L^p$ boundedness hold merely under the assumption
\eqref {log12decay}?

\end{document}